
\input amstex
\documentstyle{amsppt}

\vcorrection{-0.5cm}
\hcorrection{+0.3cm}
\pageheight{20.7cm}
\pagewidth{12.7cm}
\magnification=\magstep1
\parskip=6pt

\define\ce{{\Cal E}}
\define\cf{{\Cal F}}
\define\cg{{\Cal G}}

\define\cn{{\Cal N}}
\define\co{{\Cal O}}
\define\cs{{\Cal S}}
\define\ct{{\Cal T}}

\define\pp{{\Bbb P}}
\define\qq{{\Bbb Q}}
\define\rr{{\Bbb R}}

\define\zz{{\Bbb Z}}

\define\he{{H(\Cal E)}}

\define\hf{{H(\Cal F)}}

\define\pic{\operatorname{Pic}}
\define\rk{\operatorname{rank}}

\topmatter

\title 
  A generalization of curve genus for ample vector bundles, I
\endtitle

\author
  Hironobu Ishihara
\endauthor

\address 
  Department of Mathematics,         
  Tokyo Institute of Technology,   
  Oh-okayama, Meguro,    
  Tokyo 152, Japan               
\endaddress

\email
  ishihara\@math.titech.ac.jp
\endemail



\abstract
  A new genus $g=g(X,\ce)$ is defined for the pairs $(X,\ce)$ 
  that consist of $n$-dimensional compact complex manifolds $X$ 
  and ample vector bundles $\ce$ of rank $r$ less than $n$ on $X$.
  In case $r=n-1$, $g$ is equal to curve genus.
  Above pairs $(X,\ce)$ with $g$ less than two are classified.
  For spanned $\ce$ it is shown that $g$ is greater than or equal to
  the irregularity of $X$, and its equality condition is given.
\endabstract

\endtopmatter

\document

\subhead
  Introduction
\endsubhead

Let $X$ be a compact complex manifold of dimension $n$
and $\ce$ an ample vector bundle of rank $r$ on $X$.
A pair $(X,\ce)$ is called a generalized polarized manifold.
It is a natural and interesting problem to generalize 
sectional genus, which is defined for polarized manifolds,
for generalized polarized manifolds.
Fujita \cite{F3} defined $c_1$-sectional genus and
$\co(1)$-sectional genus for pairs $(X,\ce)$
as the sectional genus of $(X,\det\ce)$ and $(\pp_X(\ce),\he)$ 
respectively.
($\he$ is the tautological line bundle on $\pp_X(\ce)$.)
Fujita \cite{F3} also gave the classification of $(X,\ce)$ of
$c_1$-sectional genus or $\co(1)$-sectional genus less than three.

Ballico \cite{Ba} defined another sectional genus,
which was called curve genus in \cite{LMS},
for above $(X,\ce)$ with $r=n-1$.
He gave a classification of $(X,\ce)$ of curve genus zero,
and of curve genus one for spanned $\ce$ of rank two.
Recently Lanteri and Maeda \cite{LM2} have provided 
the classification of $(X,\ce)$ of curve genus zero or one
under assumption on the existence of a regular section of $\ce$.
Later Maeda \cite{M} has completed
the classification of $(X,\ce)$ of curve genus zero or one
by removing the above assumption for $\ce$.
Lanteri, Maeda and Sommese \cite{LMS} obtained that
for spanned $\ce$ curve genus is greater than or equal to 
the irregularity $q(X)$ of $X$,
and gave its equality condition. 

In the present paper we generalize curve genus
and define a new genus (which might be called ``$c_r$-sectional genus'')
$g(X,\ce)$ for above $(X,\ce)$ with $r<n$ by the formula
  $$ 2g(X,\ce)-2:=(K_X+(n-r)c_1(\ce))c_1(\ce)^{n-r-1}c_r(\ce), $$
where $K_X$ is the canonical bundle of $X$.
We note that this genus is equal to curve genus when $r=n-1$ 
and equal to ($c_1$-)sectional genus when $r=1$. 
We see that the properties for curve genus
proved in the literatures above are extended to
our $c_r$-sectional genus.

The contents of this paper are as follows.
In Section 1 we examine fundamental properties of $g(X,\ce)$.
We see that $g(X,\ce)$ is a non-negative integer and
the classifications of $(X,\ce)$ with $g(X,\ce)\le 1$ are given.
In particular we prove that if $g(X,\ce)=0$, then
$(X,\ce)\simeq(\pp^4,\co(1)^{\oplus2})$
unless $(X,\ce)$ is in two ``boundary'' cases,
namely the cases $r=1$ and $r=n-1$. 
In Section 2 we show that $g(X,\ce)\ge q(X)$ for spanned $\ce$ 
and we give its equality condition;
the result is that $(X,\ce)\simeq(\pp^4,\co(1)^{\oplus2})$
when $1<r<n-1$.
We give also another equality conditions as in \cite{LMS}.

We will further study the properties of $g(X,\ce)$ in a forthcoming paper.

\subhead
  Acknowledgment
\endsubhead

The author is grateful to Dr\. Yoshiaki Fukuma;
this paper has grown out from their conversations.

\subhead 
 1. Fundamental properties of $g(X,\ce)$
\endsubhead

Basically we use the standard notation as that in \cite{H}.
The tensor products of line bundles are usually denoted additively,
while we use multiplicative notation for intersection products.
For a vector bundle $\ce$ on a variety $X$, we denote by 
$\he$ the tautological line bundle on $\pp_X(\ce)$,
the associated projective space bundle in the sense of Grothendieck.
We say that $\ce$ is spanned if $\he$ is spanned (by global sections).
For an embedding $\iota:Y \hookrightarrow X$,
we often denote the pull-back $\iota^*\ce$ by $\ce_Y$.

\definition{Definition 1.1}
  Let $X$ be a compact complex manifold of dimension $n\ge 2$
  and $\ce$ an ample vector bundle of rank $r<n$ on $X$.
  We define a rational number $g(X,\ce)$ for a generalized
  polarized manifold $(X,\ce)$ by the formula
  $$ 2g(X,\ce)-2:=(K_X+(n-r)c_1(\ce))c_1(\ce)^{n-r-1}c_r(\ce), $$
  where $K_X$ is the canonical bundle of $X$.
\enddefinition

\example{Remark 1.2}
  Let $(X,\ce)$ be as above.
  When $r=1$, $g(X,\ce)$ is the sectional genus of a polarized manifold
  $(X,\ce)$;
  when $r=n-1$, $g(X,\ce)$ is the curve genus of a generalized
  polarized manifold $(X,\ce)$.
  In view of \cite{FulL, Theorem I}, we see that for every $r<n$,
  $$ \align
     &(K_X+(n-r)c_1(\ce))c_1(\ce)^{n-r-1}c_r(\ce)\\
  \le&(K_X+(n-1)c_1(\ce))c_1(\ce)^{n-r-1}c_r(\ce)\\
  \le&(K_X+(n-1)c_1(\ce))c_1(\ce)^{n-1}
     \endalign
  $$
  because $K_X+(n-1)\det\ce$ is nef unless $r=1$.
  Hence we obtain $g(X,\ce)\le g(X,\det\ce)$, i.e.,
  $g(X,\ce)$ is not greater than the $c_1$-sectional genus of $(X,\ce)$.
\endexample

\example{Remark 1.3}
  Let $(X,\ce)$ be as above.
  Suppose that $(X,\ce)$ satisfies the condition
  \roster
  \item"($\ast$)" there exists a section $s\in H^0(X,\ce)$ whose zero 
            locus $Z:=(s)_0$ is a smooth submanifold of $X$
            of the expected dimension $n-r$.
 \endroster
 Then $Z=c_r(\ce)$ in the Chow ring of $X$ and
  $$ \align
     &(K_X+(n-r)c_1(\ce))c_1(\ce)^{n-r-1}c_r(\ce)\\
    =&([K_X+\det\ce]_Z+(n-r-1)c_1(\ce_Z))c_1(\ce_Z)^{n-r-1}\\
    =&2g(Z,\det\ce_Z)-2.
    \endalign
  $$
  Hence we obtain $g(X,\ce)=g(Z,\det\ce_Z)$, the sectional genus
  of $(Z,\det\ce_Z)$.
  In particular, if $\ce$ is spanned, then $\ce$ satisfies ($\ast$)
  and $g(X,\ce)=g(Z,\det\ce_Z)$ because of Bertini's theorem.
\endexample

\proclaim{Proposition 1.4}
  Let $(X,\ce)$ be as in $(1.1)$.
  Then $g(X,\ce)$ is an integer.
\endproclaim

\demo{Proof}
  (Cf. \cite{Ba, Remark 4.1}.)
  We take a sufficiently ample line bundle $L\in\pic X$
  such that $\ce\otimes L$ is spanned. 
  Then there exists a non-zero section $s\in H^0(\ce\otimes L)$
  whose zero locus $Z:=(s)_0$ is a smooth submanifold of $X$
  of dimension $n-r$.
  We see that $Z=c_r(\ce\otimes L)$ in the Chow ring and
  $$ \align
     &(K_X+(n-r)c_1(\ce\otimes L))c_1(\ce\otimes L)^{n-r-1}
                                  c_r(\ce\otimes L)\\
    =&([K_X+\det(\ce\otimes L)]_Z+(n-r-1)c_1([\ce\otimes L]_Z))
      c_1([\ce\otimes L]_Z)^{n-r-1}\\
    =&2g(Z,\det[\ce\otimes L]_Z)-2\in 2\zz.
    \endalign
  $$ 
  On the other hand we have
  $$ \align
     &(K_X+(n-r)c_1(\ce\otimes L))c_1(\ce\otimes L)^{n-r-1}
                                  c_r(\ce\otimes L)\\
    =&(K_X+(n-r)(c_1(\ce)+rL))(c_1(\ce)+rL)^{n-r-1}
     (\sum_{i=0}^r c_i(\ce)c_1(L)^{r-i}).
     \endalign
  $$
  Here we can take $L=2A$ for some $A\in\pic X$, hence
  $$ \align
     &(K_X+(n-r)(c_1(\ce)+rL))(c_1(\ce)+rL)^{n-r-1}
      (\sum_{i=0}^r c_i(\ce)c_1(L)^{r-i})\\
    \equiv&(K_X+(n-r)c_1(\ce))c_1(\ce)^{n-r-1}c_r(\ce)~(\text{mod}~2)\\
    =&2g(X,\ce)-2.
    \endalign
  $$
  Thus $g(X,\ce)\in\zz$. \qed
\enddemo

\proclaim{Theorem 1.5}
  Let $X$ be a compact complex manifold of dimension $n\ge 2$
  and $\ce$ an ample vector bundle of rank $r<n$ on $X$.
  Then we have $g(X,\ce)\ge 0$. 
  Moreover, when $1<r<n-1$, we have $g(X,\ce)=0$ if and only if
            $(X,\ce)\simeq(\pp^4,\co(1)^{\oplus2})$.
\endproclaim

\demo{Proof}
  (Cf. \cite{Ba, Remark 4.2}.) 
  First we note that if $(X,\ce)\simeq(\pp^4,\co(1)^{\oplus2})$, 
  then we easily see $g(X,\ce)=0$. 
  Suppose that $g(X,\ce)\le 0$ in the following. 
  We have 
  $$ (K_X+(n-r)c_1(\ce))c_1(\ce)^{n-r-1}c_r(\ce)<0. $$
  In particular, we have $K_X+(n-r)\det\ce\not=\co_X$.
  If $K_X+(n-r)\det\ce$ is nef, by the base point free theorem,
  we get a non-zero member $D\in|m(K_X+(n-r)\det\ce)|$
  for some $m>0$.
  For every irreducible component $D'$ of $D$, we get
  $$ D'c_1(\ce)^{n-r-1}c_r(\ce)=c_1(\ce_{D'})^{n-r-1}c_r(\ce_{D'})>0 $$
  by \cite{BlG, Theorem 2.5}.
  It follows that $Dc_1(\ce)^{n-r-1}c_r(\ce)>0$, a contradiction. 
  Hence $K_X+(n-r)\det\ce$ is not nef.
  Then there exists an extremal ray $R$ such that  
  $(K_X+(n-r)\det\ce)R<0$.
  We take a rational curve $C$ belonging to $R$ 
  and satisfying $l(R)=(-K_X)C$, where $l(R)$ is the length of $R$.
  Since $l(R)\le n+1$ and
  $$ (-K_X)C>(n-r)c_1(\ce)C\ge (n-r)r\ge n-1, $$ 
  we see that $r=1$, $r=n-1$, or $r=2$ and $n=4$.

  When $r=1$, by \cite{F2, Corollary 1}, we obtain $g(X,\ce)=0$;
  when $r=n-1$, by \cite{M}, we obtain $g(X,\ce)=0$.
  When $r=2$ and $n=4$, we have $l(R)=5=n+1$, 
  hence $\pic X\simeq\zz$ and $-K_X$ is ample by 
  \cite{W, (2.4) Proposition}.
  We take the ample generator $H$ of $\pic X$ and we set 
  $\det\ce=aH$ for some $a\in\zz$.
  Since $c_1(\ce)C=2$, we have $a=1$ or 2.
  If $a=2$, then $K_X+4H$ is not nef, hence $(X,H)\simeq(\pp^4,\co(1))$
  by \cite{F2, Theorem 1}.
  Then we have $(\det\ce)L=2$ for every line $L$ in $\pp^4$.
  It follows that $\ce|_L=\co_L(1)^{\oplus2}$,
  hence we obtain $\ce=\co_{\pp^4}(1)^{\oplus2}$ 
  by \cite{OSS, Chap.~I (3.2.1)}.
  If $a=1$, we can set $-K_X=b\cdot \det\ce$ for some $b\in\zz$.
  Then we have $5=(-K_X)C=b\cdot c_1(\ce)C=2b$, a contradiction.
  This completes the proof. \qed
\enddemo

\example{Remark 1.6}
  Let $(X,\ce)$ be as in (1.5).
  When $r=1$, by Fujita \cite{F2, Corollary 1} or Ionescu 
  \cite{I, Corollary 8},
  $g(X,\ce)=0$ if and only if $(X,\ce)$ is either 
  $(\pp^n,\co(1))$, $(\pp^2,\co(2))$, $(\qq^n,\co(1))$, 
  or $(X,\ce)\simeq(\pp_{\pp^1}(\cf),\hf)$
     for some ample vector bundle $\cf$ of rank $n$ on $\pp^1$.
  ($\qq^n$ is a smooth hyperquadric in $\pp^{n+1}$.)
    
  When $r=n-1$, by Lanteri and Maeda (\cite{LM2}, \cite{M}), 
  $g(X,\ce)=0$ if and only if $(X,\ce)$ is one of the following:
  \roster 
  \item $(\pp^n,\co(1)^{\oplus(n-1)})$;
  \item $(\pp^n,\co(1)^{\oplus(n-2)}\oplus\co(2))$;
  \item $(\qq^n,\co(1)^{\oplus(n-1)})$;
  \item $X\simeq\pp_{\pp^1}(\cf)$ for some vector bundle $\cf$
                of rank $n$ on $\pp^1$, and 
                $\ce=\oplus_{j=1}^{n-1}(H(\cf)\mathbreak
                                              +\pi^*\co_{\pp^1}(b_j))$,
                where $\pi:X\to\pp^1$ is the bundle projection.
  \endroster
\endexample

\proclaim{Theorem 1.7}
  Let $X$ be a compact complex manifold of dimension $n\ge 2$
  and $\ce$ an ample vector bundle of rank $r<n$ on $X$.
  When $1<r<n-1$, we have $g(X,\ce)=1$ if and only if
            $(X,\ce)$ is one of the following: 
  \roster
  \item $(\pp^5,\co(1)^{\oplus2})$;
  \item $(\pp^5,\co(1)^{\oplus3})$;
  \item $(\qq^4,\co(1)^{\oplus2})$.
  \endroster
\endproclaim

\demo{Proof}
  The following argument is similar to that in (1.5).
  First we note that if $(X,\ce)$ is one of the above cases 
  (1),(2) and (3), then we find that $g(X,\ce)=1$ by computations.
  Suppose that $g(X,\ce)=1$ on the contrary.
  Then we have 
  $$ (K_X+(n-r)c_1(\ce))c_1(\ce)^{n-r-1}c_r(\ce)=2g(X,\ce)-2=0. $$
 
  If $K_X+(n-r)\det\ce\not=\co_X$, we infer that $K_X+(n-r)\det\ce$
  is not nef as before.
  Then we find that $r=1$, $r=n-1$, or $r=2$ and $n=4$
  by using an extremal ray.
  When $r=2$ and $n=4$, we get $g(X,\ce)\not=1$
  by the same argument as that in the proof of (1.5).
  Thus it remains to consider the case $K_X+(n-r)\det\ce=\co_X$.

  Suppose that $K_X+(n-r)\det\ce=\co_X$ and $1<r<n-1$ in the following.
  Since $K_X+(n-r-1)\det\ce$ is not nef,
  there exists an extremal rational curve $C$ such that
  $(K_X+(n-r-1)\det\ce)C<0$ and $l(R)=(-K_X)C$, where $R:=\rr_+[C]$.
  Since 
  $$ n+1\ge (-K_X)C=(n-r)c_1(\ce)C\ge (n-r)r\ge 2(n-2), $$ 
  we see that $(n,r)=(4,2),(5,2)$, or $(5,3)$.

  When $(n,r)=(5,2)$, we have $l(R)=6=n+1$, 
  hence $\pic X\simeq\zz$ and $-K_X$ is ample.
  We set $\det\ce=aH$ for the ample generator $H$ of $\pic X$
  and $a\in\zz$;
  then $a=1$ or 2 since $c_1(\ce)C=2$.
  In case $a=2$, we have $K_X+6H=0$, hence $(X,H)\simeq(\pp^5,\co(1))$
  by the Kobayashi-Ochiai theorem \cite{KO}.
  Then we get $\ce\simeq\co(1)^{\oplus2}$ since $\det\ce=\co(2)$. 
  This is the case (1) of our theorem.
  In case $a=1$, we have $\pic X=\zz\cdot\det\ce$ and
  $-K_X=3\det\ce$.
  
  Let $p:\pp_X(\ce)\to X$ be the projective bundle associated to $\ce$.
  For $P:=\pp_X(\ce)$, we see that $-K_P=2\he+p^*(2\det\ce)$ is ample
  and $P$ has two extremal rays.
  One of them is generated by the class of a fiber of $p$;
  let $R'$ be another extremal ray and 
  $f_{R'}:P\to W$ the contraction morphism of $R'$. 
  We take an extremal rational curve $C'$ belonging to $R'$
  and satisfying $l(R')=(-K_P)C'$.
  We have
  $$ (-K_P)C'=2\he C'+p^*(2\det\ce)\cdot C'\ge 6 $$
  since $p(C')$ is a rational curve on $X$.
  Hence from \cite{W, (2.4) Proposition}, we infer that $l(R')=6$, 
  $W$ is a curve, $\pic F\simeq\zz$ and $-K_F$ is ample 
  for a general fiber $F$ of $f_{R'}$.
  It follows that $\he C'=1$ and $p^*(2\det\ce)\cdot C'=4$, 
  hence we get 
  $[p^*(2\det\ce)]_F=4\he|_F$ and $-K_F=(-K_P)|_F=6\he|_F$.
  Then we get $(F,\he|_F)\simeq(\pp^5,\co(1))$ 
  by the Kobayashi-Ochiai theorem.
  Since $p(C')$ is a curve on $X$, we infer that
  $p|_F:F\to X$ is finite and surjective.
  Hence we get $X\simeq\pp^5$ by \cite{Laz, Theorem 4.1}.
  This is a contradiction since $\pic X\simeq\zz\cdot\det\ce$.
  
  When $(n,r)=(4,2)$ or $(5,3)$, we have $r=n-2$, 
  $K_X+\det\ce$ is not nef, and $(-K_X-\det\ce)C=(\det\ce)C\ge 2$.
  Hence by \cite{Z, Proposition 1.1'}, $(X,\ce)$ is one of the following:
  \roster
  \item"(i)" $(\pp^n,\co(1)^{\oplus(n-2)})$;
  \item"(ii)" $(\pp^n,\co(1)^{\oplus(n-3)}\oplus\co(2))$;
  \item"(iii)" $(\qq^n,\co(1)^{\oplus(n-2)})$;
  \item"(iv)" $X\simeq\pp_B(\cf)$ and $\ce=\hf\otimes\pi^*\cg$,
              where $\cf$ and $\cg$ are vector bundles
              on a smooth curve $B$ such that $\rk\cf=n,\rk\cg=n-2$,
              and $\pi$ is the bundle projection $X\to B$.
  \endroster
  In case (i) (resp.~(iii)), we get $n=5$ (resp~$n=4$) from 
  $K_X+2\det\ce=\co_X$, which leads to the case (2) (resp.~(3))
  of our theorem.
  Case (ii) is ruled out by $K_X+2\det\ce=\co_X$ and $n\ge 4$.
               
  In case (iv), we get $n=4$ and $K_B+\det\cf+2\det\cg=\co_B$
  from $K_X+2\det\ce=\co_X$.
  Since $\det\ce=2\hf+\pi^*\det\cg$ is ample, we have 
  $$ 0<(\det\ce)^4=16\hf^4+32\hf^3(\pi^*\det\cg)=16c_1(\cf)+32c_1(\cg). 
  $$
  Hence we find that $\deg K_B<0$, and then $B\simeq\pp^1$ and 
  $c_1(\cf)+2c_1(\cg)=2$.
  Since $\cf$ and $\cg$ are vector bundles on $\pp^1$, we can write
  $\cf=\oplus_{i=1}^4\co_{\pp^1}(a_i)$ and
  $\cg=\co_{\pp^1}(b_1)\oplus\co_{\pp^1}(b_2)$ for some 
  $a_i,b_1,b_2\in\zz,a_1\le\dots\le a_4$.
  A natural surjection $\cf\to\co_{\pp^1}(a_1)$ determines 
  a section $Z:=\pp(\co_{\pp^1}(a_1))$ of $\pi$.
  Since $\hf|_Z=\co_{\pp^1}(a_1)$, we have
  $\ce_Z=\co_{\pp^1}(a_1+b_1)\oplus\co_{\pp^1}(a_1+b_2)$,
  hence $a_1+b_1>0$ and $a_1+b_2>0$.
  It follows that 
  $$ 2=c_1(\cf)+2c_1(\cg)=\sum_{i=1}^4 a_i+2b_1+2b_2\ge 4, $$
  a contradiction.
  We have thus proved the theorem. \qed
\enddemo

\example{Remark 1.8}
  Let $(X,\ce)$ be as in (1.7). 
  When $r=1$, by Fujita \cite{F2, Corollary 2} or Ionescu
  \cite{I, Corollary 9}, 
  $g(X,\ce)=1$ if and only if
  $(X,\ce)$ is either a Del Pezzo manifold (i.e. $K_X+(n-1)\ce=\co_X$)
  or $(X,\ce)\simeq(\pp_B(\cf),\hf)$ for some ample vector bundle 
  $\cf$ of rank $n$ on an elliptic curve $B$. 
  Note that Del Pezzo manifolds have been classified in \cite{F1}.
 
  When $r=n-1$, Lanteri and Maeda (\cite{LM2}, \cite{M}) have shown 
  that if $g(X,\ce)=1$, then $(X,\ce)$ is one of the following:
  \roster
  \item $(\pp^n,\co(1)^{\oplus(n-2)}\oplus\co(3))$;
  \item $(\pp^n,\co(1)^{\oplus(n-3)}\oplus\co(2)^{\oplus2})$;
  \item $(\pp^3,\cn(2))$, where $\cn$ is the null correlation 
                    bundle on $\pp^3$ (see \cite{OSS, p.~76});
  \item $(\qq^n,\co(1)^{\oplus(n-2)}\oplus\co(2))$;
  \item $(\qq^3,\cs(2))$, where $\cs$ is the spinor bundle 
                  on $\qq^3$ (see \cite{Ot, Definition 1.3});     
  \item $(\qq^4,\cs(2)\oplus\co(1))$;
  \item $(\pp^2\times\pp^1,\co(1,1)\oplus\co(2,1))$;
  \item $(\pp^2\times\pp^1,p_1^*\ct_{\pp^2}\otimes\co(0,1))$,
                 where $\ct_{\pp^2}$ is the tangent bundle of $\pp^2$
                 and $p_1$ is the first projection 
                 $\pp^2\times\pp^1\to\pp^2$;
  \item $X$ is a Fano manifold of index $n-1$ with $\pic X\simeq\zz$, 
                  and $\ce=H^{\oplus(n-1)}$ for the ample generator $H$ 
                  of $\pic X$;
  \item $X\simeq\pp_B(\cf)$ for some vector bundle $\cf$
                  of rank $n$ on an elliptic curve $B$, and
                  $\ce_F=\co_{\pp^{n-1}}(1)^{\oplus(n-1)}$ 
                  for every fiber $F$ of the bundle projection $X\to B$;
  \item $n=4$ and $\pp_X(\ce)$ has two projective $\pp^2$-bundle 
                   structures over smooth Fano 4-folds 
                   of index 1 with $b_2=1$ and pseudoindex $\ge 3$.
  \endroster
\endexample

\subhead
  2. A relation between $g(X,\ce)$ and $q(X)$
\endsubhead

\proclaim{Proposition 2.1}
  Let $X$ be a compact complex manifold of dimension $n\ge 2$
  and $A$ a spanned ample line bundle on $X$.
  Then we have $g(X,A)\ge q(X)$, and equality holds if and only if
  $(X,A)$ is one of the following:
  \roster
  \item $(\pp^n,\co(1))$;
  \item $(\pp^2,\co(2))$;
  \item $(\qq^n,\co(1))$;
  \item $(\pp_B(\cf),\hf)$, where $\cf$ is a spanned
        ample vector bundle of rank $n$ on a smooth curve $B$.
  \endroster
\endproclaim

\demo{Proof}
  See, e.g., \cite{BeS, Theorem 7.2.10}. \qed
\enddemo

\proclaim{Theorem 2.2}
  Let $X$ be a compact complex manifold of dimension $n\ge 2$
  and $\ce$ a spanned ample vector bundle of rank $r<n$ on $X$.
  Then we have $g(X,\ce)\ge q(X)$. 
  Moreover, when $1<r<n-1$, we have $g(X,\ce)=q(X)$ if and only if
  $(X,\ce)\simeq(\pp^4,\co(1)^{\oplus2})$.
\endproclaim

\demo{Proof}
  (Cf. \cite{LMS, Introduction}.) 
  Since $\ce$ is spanned, by (1.3), there exists a nonzero section
  $s\in H^0(X,\ce)$ whose zero locus $Z:=(s)_0$ is a smooth submanifold
  of $X$ of dimension $n-r$.
  Then a natural map $H^1(X,\zz)\to H^1(Z,\zz)$ is injective
  by a Lefschetz-type theorem \cite{LM1, Theorem 1.3}.
  It follows that $H^1(X,\co_X)\to H^1(Z,\co_Z)$ is also injective
  and $q(X)\le q(Z)$.
  On the other hand, we have $g(X,\ce)=g(Z,\det\ce_Z)$ and
  $\det\ce_Z$ is ample and spanned.
  Hence we get $g(Z,\det\ce_Z)\ge q(Z)$ by (2.1), 
  thus $g(X,\ce)\ge q(X)$.

  Suppose that $g(X,\ce)=q(X)$ and $1<r<n-1$ in the following.
  Then we get $g(Z,\det\ce_Z)=q(Z)$ by the argument above.
  Since $\det\ce_Z$ is ample and spanned,
  $(Z,\det\ce_Z)$ is one of the four cases in (2.1).
  If $(Z,\det\ce_Z)$ is the case (4) in (2.1), then
  $(\det\ce_Z)C=1$ for some curve $C$ in a fiber of the bundle map
  $Z\to B$, which is impossible.
  If $(Z,\det\ce_Z)$ is one of the other cases in (2.1),
  we get $q(Z)=0$ and then $g(X,\ce)=0$.
  It follows that $(X,\ce)\simeq(\pp^4,\co(1)^{\oplus2})$ by (1.5).
  On the other hand, it is easy to see that
  $g(\pp^4,\co(1)^{\oplus2})=0=q(\pp^4)$.
  We have thus proved the theorem. \qed
\enddemo

\example{Remark 2.3}
  Let $(X,\ce)$ be as in (2.2).
  When $r=n-1$, Lanteri, Maeda and Sommese 
  \cite{LMS, Theorem} have shown that
  $g(X,\ce)=q(X)$ if and only if $(X,\ce)$ is one of the following:
  \roster
  \item $(\pp^n,\co(1)^{\oplus(n-1)})$;
  \item $(\pp^n,\co(1)^{\oplus(n-2)}\oplus\co(2))$;
  \item $(\qq^n,\co(1)^{\oplus(n-1)})$;
  \item $X\simeq\pp_B(\cf)$ and $\ce=\hf\otimes\pi^*\cg$,
        where $\cf$ and $\cg$ are vector bundles on a smooth curve $B$
        with $\rk\cf=n$ and $\rk\cg=n-1$, and $\pi:X\to B$ is the 
        bundle projection.
  \endroster
\endexample

\proclaim{Corollary 2.4}
  Let $(X,\ce)$ be as in {\rm(2.2)}.
  Then the following conditions are equivalent.
  \roster
  \item"(i)" $g(X,\ce)=q(X)$.
  \item"(ii)" $h^0(m(K_X+(n-r)\det\ce))=0$ for all $m\ge 1$.
  \item"(iii)" $K_X+(n-r)\det\ce$ is not nef.
  \endroster
\endproclaim

\demo{Proof}
  The assertion holds when $r=1$ (resp.~$r=n-1$) 
  by \cite{S, (4.1) and (4.2)} (resp.~\cite{LMS, (2.1) Corollary}).
  Hence we may assume that $1<r<n-1$.
  If (i) holds, then we easily find that (ii) holds by (2.2).
  If (ii) holds, then (iii) holds by the base point free theorem.
  If (iii) holds, then we find that $g(X,\ce)=0=q(X)$
  by the argument in (1.5). \qed
\enddemo

\enddocument